\documentclass[12pt,leqno]{article}

\def\dist{\mathop{\rm dist}}
\def\Lip{\mathop{\rm Lip}}

\newtheorem{theorem}{Theorem}
\newtheorem{lemma}[theorem]{Lemma}
\newtheorem{proposition}[theorem]{Proposition}
\newtheorem{sublemma}[theorem]{Sublemma}
\newtheorem{definition}[theorem]{Definition}
\newtheorem{corollary}[theorem]{Corollary}
\newtheorem{problem}[theorem]{Problem}
\newtheorem{remark}[theorem]{Remark}
\newtheorem{claim}[theorem]{Claim}
\newtheorem{assumptions}[theorem]{Assumptions}
\newtheorem{examples}[theorem]{Examples}
\newtheorem{basicfact}[theorem]{Basic Fact}

\newcommand{\begintheorem}{\addtocounter{equation}{1}\begin{theorem}}
\newcommand{\beginlemma}{\addtocounter{equation}{1}\begin{lemma}}
\newcommand{\beginproposition}{\addtocounter{equation}{1}\begin{proposition}}
\newcommand{\beginsublemma}{\addtocounter{equation}{1}\begin{sublemma}}
\newcommand{\begindefinition}{\addtocounter{equation}{1}\begin{definition}}
\newcommand{\begincorollary}{\addtocounter{equation}{1}\begin{corollary}}
\newcommand{\beginproblem}{\addtocounter{equation}{1}\begin{problem}}
\newcommand{\beginremark}{\addtocounter{equation}{1}\begin{remark}}
\newcommand{\beginclaim}{\addtocounter{equation}{1}\begin{claim}}
\newcommand{\beginassumptions}{\addtocounter{equation}{1}\begin{assumptions}}
\newcommand{\beginexamples}{\addtocounter{equation}{1}\begin{examples}}
\newcommand{\beginbasicfact}{\addtocounter{equation}{1}\begin{basicfact}}

\begin{document}

\title{A few aspects of analysis on metric spaces}


\author{Stephen Semmes}

\date{}

\maketitle




%

\tableofcontents

\bigskip\bigskip

\section{Metric spaces and Lipschitz functions}
\label{Metric spaces and Lipschitz functions}
\setcounter{equation}{0}

	Let $(M, d(x,y))$ be a metric space.  Thus $M$ is a nonempty set,
$d(x,y)$ is a nonnegative real-valued function on $M \times M$ which
is equal to $0$ exactly when $x = y$, $d(x,y)$ is symmetric in $x$ and
$y$, so that $d(x,y) = d(y,x)$ for all $x, y \in M$, and $d(x,y)$
satisfies the triangle inequality,
\begin{equation}
\label{triangle inequality}
	d(x,z) \le d(x,y) + d(y,z)
\end{equation}
for all $x$, $y$, and $z$ in $M$.  If $x$ is an element of $M$ and $r$
is a positive real number, then we write $B(x,r)$ and
$\overline{B}(x,r)$ for the open and closed balls in $M$ with center
$x$ and radius $r$, i.e.,
\begin{equation}
	B(x,r) = \{y \in M : d(x,y) < r\}, \ 
		\overline{B}(x,r) = \{y \in M : d(x,y) \le r \}.
\end{equation}

	Suppose that $f(x)$ is a real or complex-valued function on
$M$, and that $L$ is a nonnegative real number.  We say that $f$ is
\emph{$L$-Lipschitz} if
\begin{equation}
\label{def of L-Lipschitz}
	|f(x) - f(y)| \le L \, d(x,y)
\end{equation}
for all $x, y \in M$.  We also simply say that $f$ is Lipschitz if
it is $L$-Lipschitz for some $L$.  If $f$ is Lipschitz, then we define
$\|f\|_{\Lip }$ to be the supremum of
\begin{equation}
	\frac{|f(x) - f(y)|}{d(x,y)}
\end{equation}
over all $x, y \in M$, where this ratio is replaced with $0$ when $x =
y$.  Thus $f$ is $\|f\|_{\Lip }$-Lipschitz when $f$ is Lipschitz, and
this is the smallest choice of $L$ for which $f$ is $L$-Lipschitz.
Note that $\|\cdot \|_{\Lip }$ is a seminorm, so that
\begin{equation}
	\|a \, f + b \, g \|_{\Lip } 
		\le |a| \, \|f\|_{\Lip } + |b| \, \|g\|_{\Lip }
\end{equation}
for all constants $a$, $b$ and Lipschitz functions $f$, $g$ on $M$.
Also, $\|f\|_{\Lip } = 0$ if and only if $f$ is a constant function on $M$.

	If $f$ and $g$ are real-valued $L$-Lipschitz functions on $M$,
then $\max(f,g)$ and $\min(f,g)$ are $L$-Lipschitz functions too.  Let
us check this for $\max(f,g)$.  It is enough to show that
\begin{equation}
	\max(f,g)(x) - \max(f,g)(y) \le L \, d(x,y)
\end{equation}
for all $x, y \in M$, i.e., one can interchange the roles of $x$ and
$y$ to get an inequality in the opposite direction.  Assume, for the
sake of definiteness, that $\max(f,g)(x) = f(x)$.  Then we have
\begin{equation}
	\ \max(f,g)(x) = f(x) \le f(y) + L \, d(x,y) 
			\le \max(f,g)(y) + L \, d(x,y),
\end{equation}
which is what we wanted.

	Here is a generalization of this fact.

\beginlemma
\label{lemma about sup of a family of L-Lipschitz functions}
Let $\{f_\sigma\}_{\sigma \in A}$ be a family of real-valued functions
on $M$ which are all $L$-Lipschitz for some $L \ge 0$.  Assume also
that there is point $p$ in $M$ such that the set of real numbers
$\{f_\sigma(p) : \sigma \in A\}$ is bounded from above.  Then the set
$\{f_\sigma(x) : \sigma \in A\}$ is bounded from above for every $x$
in $M$ (but not uniformly in $x$ in general), and $\sup \{f_\sigma(x)
: \sigma \in A\}$ is an $L$-Lipschitz function on $M$.
\end{lemma}

	Indeed, because $f_\sigma$ is $L$-Lipschitz for all $\sigma$
in $A$, we have that
\begin{equation}
	f_\sigma(x) \le f_\sigma(y) + L \, d(x,y)
\end{equation}
for all $x$, $y$ in $M$.  Applying this to $y = p$, we see that 
$\{f_\sigma(x) : \sigma \in A\}$ is bounded from above for every
$x$, because of the corresponding property for $p$.  If $F(x) =
\sup \{f_\sigma(x) : \sigma \in A\}$, then
\begin{equation}
	F(x) \le F(y) + L \, d(x,y)
\end{equation}
for all $x$, $y$ in $M$, so that $F$ is $L$-Lipschitz on $M$.  

	For the record, let us write down the analogous statement
for infima of $L$-Lipschitz functions.

\beginlemma
\label{lemma about inf of a family of L-Lipschitz functions}
Let $\{f_\sigma\}_{\sigma \in A}$ be a family of real-valued functions
on $M$ which are all $L$-Lipschitz for some $L \ge 0$.  Assume also
that there is point $q$ in $M$ such that the set of real numbers
$\{f_\sigma(q) : \sigma \in A\}$ is bounded from below.  Then the set
$\{f_\sigma(x) : \sigma \in A\}$ is bounded from below for every $x$
in $M$, and $\inf \{f_\sigma(x) : \sigma \in A\}$ is an $L$-Lipschitz
function on $M$.
\end{lemma}

	For any point $w$ in $M$, $d(x,w)$ defines a $1$-Lipschitz
function of $x$ on $M$.  This can be shown using the triangle inequality.
Suppose now that $f(x)$ is an $L$-Lipschitz function on $M$.  For
each $w \in M$, define $f_w(x) = f(w) + L \, d(x,w)$.  The fact
that $f$ is $L$-Lipschitz implies that
\begin{equation}
	f(x) \le f_w(x) \quad\hbox{for all } x, w \in M.
\end{equation}
Of course $f_x(x) = x$, and hence
\begin{equation}
	f(x) = \inf \{f_w(x) : w \in M\}.
\end{equation}
Each function $f_w(x)$ is $L$-Lipschitz in $x$, since $d(x,w)$ is
$1$-Lipschitz in $w$.  

	Similarly, we can set $\widetilde{f}_w(x) = f(x) - L \, d(x,w)$,
and then we have that
\begin{equation}
	f(x) = \sup \{\widetilde{f}_w(x) : w \in M \},
\end{equation}
and that $\widetilde{f}_w(x)$ is an $L$-Lipschitz function of $x$ for
every $w$.

	Here is a variant of these themes.  Let $E$ be a nonempty
subset of $M$, and suppose that $f$ is a real-valued function on $E$
which is $L$-Lipschitz, so that
\begin{equation}
	|f(x) - f(y)| \le L \, d(x,y)
\end{equation}
for all $x$, $y$ in $M$.  For each $w$ in $E$, set $f_w(x) = f(x) + L
\, d(x,w)$ and $\widetilde{f}_w(x) = f(x) - L \, d(x,w)$.  Consider
\begin{equation}
	F(x) = \inf \{f_w(x) : w \in E \}, \quad
		\widetilde{F}(x) = \sup \{\widetilde{f}_w(x) : w \in E \},
\end{equation}
for $x$ in $M$.  For the same reasons as before, $F(x) =
\widetilde{F}(x) = f(x)$ when $x$ lies in $E$.  Using Lemmas
\ref{lemma about sup of a family of L-Lipschitz functions} and
\ref{lemma about inf of a family of L-Lipschitz functions},
one can check that $F$ and $\widetilde{F}$ are $L$-Lipschitz
real-valued functions on all of $M$, i.e., they are extensions
of $f$ from $E$ to $M$ with the same Lipschitz constant $L$.

	If $H(x)$ is any other real-valued function on $M$ which
agrees with $f$ on $E$ and is $L$-Lipschitz, then
\begin{equation}
	\widetilde{f}_w(x) \le H(x) \le f_w(x)
\end{equation}
for all $w$ in $E$ and $x$ in $M$, and hence
\begin{equation}
	\widetilde{F}(x) \le H(x) \le F(x)
\end{equation}
for all $x$ in $M$.

\beginremark
\label{def of dist(x, S), and it is 1-Lipschitz}
{\rm 
If $S$ is any nonempty subset of $M$, define $\dist(x, S)$ for
$x$ in $M$ by
\begin{equation}
	\dist(x, S) = \inf_{y \in S} d(x,y).
\end{equation}
This function is always $1$-Lipschitz in $x$, by Lemma \ref{lemma
about inf of a family of L-Lipschitz functions}.
}
\end{remark}

\section{Lipschitz functions of order $\alpha$}
\label{Lipschitz functions of order alpha}
\setcounter{equation}{0}

	Let $(M, d(x,y))$ be a metric space, and let $\alpha$ be a
positive real number.  A real or complex-valued function $f$ on $M$
is said to be \emph{Lipschitz of order $\alpha$} if there is nonnegative
real number $L$ such that
\begin{equation}
\label{def of Lip of order alpha}
	|f(x) - f(y)| \le L \, d(x,y)^\alpha
\end{equation}
for all $x, y \in M$.  This reduces to the Lipschitz condition
discussed in Section \ref{Metric spaces and Lipschitz functions} when
$\alpha = 1$.  We shall sometimes write $\Lip \alpha$ for the
collection of Lipschitz functions of order $\alpha$, which might be
real or complex valued, depending on the context.  One also sometimes
refers to these functions as being ``H\"older continuous of order
$\alpha$''.

	If $f$ is Lipschitz of order $\alpha$, then we define
$\|f\|_{\Lip \alpha}$ to be the supremum of
\begin{equation}
	\frac{|f(x) - f(y)|}{d(x,y)^\alpha}
\end{equation}
over all $x, y \in M$, where this quantity is replaced with $0$ when
$x = y$.  In other words, $\|f\|_{\Lip \alpha}$ is the smallest choice
of $L$ so that (\ref{def of Lip of order alpha}) holds for all $x, y
\in M$.  This defines a seminorm on the space of Lipschitz functions
of order $\alpha$, as before, with $\|f\|_{\Lip \alpha} = 0$ if and
only if $f$ is constant.  Of course $\|f\|_{\Lip 1}$ is the same as
$\|f\|_{\Lip }$ from Section \ref{Metric spaces and Lipschitz functions}.

	If $f$ and $g$ are real-valued functions on $M$ which are
Lipschitz of order $\alpha$ with constant $L$, then $\max(f,g)$ and
$\min(f,g)$ are also Lipschitz of order $\alpha$ with constant $L$.
This can be shown in the same manner as for $\alpha = 1$.  Similarly,
the analogues of Lemmas \ref{lemma about sup of a family of
L-Lipschitz functions} and \ref{lemma about inf of a family of
L-Lipschitz functions} for Lipschitz functions of order $\alpha$
hold for essentially the same reasons as before.

	However, if $\alpha > 1$, it may be that the only functions
that are Lipschitz of order $\alpha$ are the constant functions.  This
is the case when $M = {\bf R}^n$, for instance, equipped with the
standard Euclidean metric, because a function in $\Lip \alpha$ with
$\alpha > 1$ has first derivatives equal to $0$ everywhere.  Instead
of using derivatives, it is not hard to show that the function has to
be constant through more direct calculation too.

	This problem does not occur when $\alpha < 1$.

\beginlemma
\label{if 0 < alpha le 1, then (a+b)^alpha le a^alpha + b^alpha}
If $0 < \alpha \le 1$ and $a$, $b$ are nonnegative real numbers, then
$(a+b)^\alpha \le a^\alpha + b^\alpha$.
\end{lemma}

	To see this, observe that 
\begin{equation}
	\max(a,b) \le (a^\alpha + b^\alpha)^{1/\alpha},
\end{equation}
and hence
\begin{eqnarray}
	a + b & \le & \max(a,b)^{1-\alpha} \, (a^\alpha + b^\alpha)	\\
		& \le & (a^\alpha + b^\alpha)^{1 + (1-\alpha)/\alpha}
		= (a^\alpha + b^\alpha)^{1/\alpha}.		\nonumber
\end{eqnarray}

\begincorollary
\label{d(x,y)^alpha a metric is d(x,y) is, 0 < alpha le 1}
If $(M, d(x,y))$ is a metric space and $\alpha$ is a real number
such that $0 < \alpha \le 1$, then $d(x,y)^\alpha$ also defines
a metric on $M$.
\end{corollary}

	This is easy to check.  The main point is that $d(x,y)^\alpha$
satisfies the triangle inequality, because of Lemma \ref{if 0 < alpha
le 1, then (a+b)^alpha le a^alpha + b^alpha} and the triangle
inequality for $d(x,y)$.

	A function $f$ on $M$ is Lipschitz of order $\alpha$ with
respect to the original metric $d(x,y)$ if and only if it is Lipschitz
of order $1$ with respect to $d(x,y)^\alpha$, and with the same norm.
In particular, for each $w$ in $M$, $d(x,w)^\alpha$ satisfies
(\ref{def of Lip of order alpha}) with $L = 1$ when $0 < \alpha \le
1$, because of the triangle inequality for $d(u,v)^\alpha$.

\section{Some functions on the real line}
\label{Some functions on the real line}
\setcounter{equation}{0}

	Fix $\alpha$, $0 < \alpha \le 1$.  For each nonnegative
integer $n$, consider the function
\begin{equation}
\label{2^{-n alpha} exp(2^n i x)}
	2^{-n \alpha} \exp(2^n \, i \, x)
\end{equation}
on the real line ${\bf R}$, where $\exp u$ denotes the usual exponential
$e^u$.  Let us estimate the $\Lip \alpha$ norm of this function.

	Recall that
\begin{equation}
\label{|exp(i u) - exp(i v)| le |u-v|}
	|\exp(i \, u) - \exp(i \, v)| \le |u-v|
\end{equation}
for all $u, v \in {\bf R}$.  Indeed, one can write $\exp(i \, u) -
\exp(i \, v)$ as the integral between $u$ and $v$ of the derivative of
$\exp(i \, t)$, and this derivative is $i \, \exp(i \, t)$, which has
modulus equal to $1$ at every point.

	Thus, for any $x, y \in {\bf R}$, we have that
\begin{equation}
  |2^{-n \alpha} \exp(2^n \, i \, x) - 2^{-n \alpha} \exp(2^n \, i \, y)|
	\le 2^{n (1-\alpha)} \, |x-y|.
\end{equation}
Of course
\begin{eqnarray}
\lefteqn{|2^{-n \alpha} \exp(2^n \, i \, x) 
   - 2^{-n \alpha} \exp(2^n \, i \, y)|} 	\\
	& & \le 
  2^{-n \alpha} |\exp(2^n \, i \, x)| + 2^{-n \alpha} |\exp(2^n \, i \, y)|
	= 2^{-n \alpha + 1}			\nonumber
\end{eqnarray}
as well.  As a result,
\begin{eqnarray}
\lefteqn{|2^{-n \alpha} \exp(2^n \, i \, x) 
   - 2^{-n \alpha} \exp(2^n \, i \, y)|} 	\\
	& & \le \Bigl(2^{n (1-\alpha)} \, |x-y| \Bigr)^{\alpha}
		\, \Bigl(2^{-n \alpha + 1} \Bigr)^{1-\alpha}
		= 2^{1-\alpha} \, |x-y|^\alpha.		\nonumber
\end{eqnarray}
This shows that the function (\ref{2^{-n alpha} exp(2^n i x)})
has $\Lip \alpha$ norm (with respect to the standard Euclidean
metric on ${\bf R}$) which is at most $2^{1-\alpha}$.  In the
opposite direction, if $2^n (x-y) = \pi$, then
\begin{eqnarray}
\lefteqn{|2^{-n \alpha} \exp(2^n \, i \, x) 
    - 2^{-n \alpha} \exp(2^n \, i \, y)|}		\\
	& & = 					
 2^{-n \alpha} |\exp(2^n \, i \, x)| + 2^{-n \alpha} |\exp(2^n \, i \, y)|
							\nonumber \\
	& & = 2^{-n \alpha + 1} = 2 \pi^{-\alpha} \, |x-y|^\alpha,
							\nonumber
\end{eqnarray}
so that the $\Lip \alpha$ norm is at least $2 \pi^{-\alpha}$.

	Now suppose that $f(x)$ is a complex-valued function on ${\bf
R}$ of the form
\begin{equation}
\label{f(x) = sum_{n=0}^infty a_n 2^{-n alpha} exp(2^n i x)}
	f(x) = \sum_{n=0}^\infty a_n \, 2^{-n \alpha} \exp(2^n \, i \, x),
\end{equation}
where the $a_n$'s are complex numbers.  We assume that the $a_n$'s are
bounded, which implies that the series defining $f(x)$ converges
absolutely for each $x$.  Set 
\begin{equation}
	A = \sup_{n \ge 0} |a_n|.
\end{equation}

	Let $m$ be a nonnegative integer.  For each $x$ in ${\bf R}$
we have that
\begin{equation}
\label{bound for sum from m to infty}
	\Bigl|\sum_{n=m}^\infty a_n \, 2^{-n \alpha} \exp(2^n \, i \, x) \Bigr|
	\le \sum_{n=m}^\infty A \, 2^{-n \alpha}
	= A \, (1 - 2^{-\alpha})^{-1} \, 2^{-m \alpha}.
\end{equation}
If $m \ge 1$ and $x, y \in {\bf R}$, then (\ref{|exp(i u) - exp(i v)|
le |u-v|}) yields
\begin{eqnarray}
\label{bound for difference of sums from 0 to m-1}
\lefteqn{\Bigl|\sum_{n=0}^{m-1} a_n \, 2^{-n \alpha} \exp(2^n \, i \, x) 
	- \sum_{n=0}^{m-1} a_n \, 2^{-n \alpha} \exp(2^n \, i \, y) \Bigr|} \\
   & & \le \sum_{n=0}^{m-1} A \, 2^{n (1-\alpha)} \, |x-y|
 \le A \, 2^{(m-1) (1-\alpha)} \Bigl(\sum_{j=0}^\infty 2^{-j (1-\alpha)}\Bigr)
		\, |x-y|				\nonumber \\	
   & & = A \, 2^{(m-1) (1-\alpha)} \, (1 - 2^{-(1-\alpha)})^{-1} \, |x-y|.
								\nonumber
\end{eqnarray}
Here we should assume that $\alpha < 1$, to get the convergence of
$\sum_{j=0}^\infty 2^{-j (1-\alpha)}$.  

	Fix $x, y \in {\bf R}$.  If $|x-y| > 1/2$, then we apply
(\ref{bound for sum from m to infty}) with $m = 0$ to both $x$ and $y$
to get that
\begin{eqnarray}
	|f(x) - f(y)| \le |f(x)| + |f(y)| 
		& \le & 2 \, A \, (1 - 2^{-\alpha})^{-1}		\\
	   & \le & 2^{1 + \alpha} \, A \, (1 - 2^{-\alpha})^{-1} |x-y|^\alpha.
								\nonumber
\end{eqnarray}
Assume now that $|x-y| \le 1/2$, and choose $m \in {\bf Z}_+$ so that
\begin{equation}
	2^{-m-1} < |x-y| \le 2^{-m}.
\end{equation}
Combining (\ref{bound for sum from m to infty}) and (\ref{bound for
difference of sums from 0 to m-1}), with (\ref{bound for sum from m to
infty}) applied to both $x$ and $y$, we obtain that
\begin{eqnarray}
\lefteqn{\quad |f(x) - f(y)|}						\\
  & & \le 2 \, A \, (1 - 2^{-\alpha})^{-1} \, 2^{-m \alpha}
	   + A \, 2^{(m-1) (1-\alpha)} \, (1 - 2^{-(1-\alpha)})^{-1} \, |x-y|
								\nonumber \\
  & & \le 2^{1+\alpha} \, A \, (1 - 2^{-\alpha})^{-1} \, |x-y|^\alpha
	+ A \, 2^{-(1-\alpha)} \, (1 - 2^{-(1-\alpha)})^{-1} \, |x-y|^\alpha.
								\nonumber
\end{eqnarray}
Therefore, for all $x, y \in {\bf R}$, we have that
\begin{eqnarray}
\lefteqn{|f(x) - f(y)|}						\\
  & & \le A (2^{1+\alpha} \, (1 - 2^{-\alpha})^{-1}
  + 2^{-(1-\alpha)} \, (1 - 2^{-(1-\alpha)})^{-1}) \, |x-y|^\alpha
								\nonumber
\end{eqnarray}
when $0 < \alpha < 1$.  In other words, $f$ is Lipschitz of order $\alpha$,
and
\begin{equation}
	\enspace \|f\|_{\Lip \alpha} 
		\le \Bigl(\sup_{n \ge 0} |a_n| \Bigr) \,
			(2^{1+\alpha} \, (1 - 2^{-\alpha})^{-1}
  			+ 2^{-(1-\alpha)} \, (1 - 2^{-(1-\alpha)})^{-1}).
\end{equation}

	To get an inequality going in the other direction we shall compute
as follows.  Let $\psi(x)$ be a function on ${\bf R}$ such that the
Fourier transform $\widehat{\psi}(\xi)$ of $\psi$,
\begin{equation}
   \widehat{\psi}(\xi) = \int_{\bf R} \exp (i \, \xi \, x ) \, \psi(x) \, dx
\end{equation}
is a smooth function which satisfies $\widehat{\psi}(1) = 1$ and
$\widehat{\psi}(\xi) = 0$ when $0 \le \xi \le 1/2$ and when $\xi \ge
2$.  One can do this with $\psi(x)$ in the Schwartz class of smooth
functions such that $\psi(x)$ and all of its derivatives are bounded
by constant multiples of $(1+|x|)^{-k}$ for every positive integer
$k$.

	For each nonnegative integer $j$, let us write $\psi_{2^j}(x)$
for the function $2^j \, \psi(2^j \, x)$.  Thus
\begin{equation}
	\widehat{\psi_{2^j}}(\xi) = \widehat{\psi}(2^{-j} \, \xi).
\end{equation}
In particular, $\widehat{\psi_{2^j}}(2^j) = 1$, and $\widehat{\psi_{2^j}}(2^l)
= 0$ when $l$ is a nonnegative integer different from $j$.  Hence
\begin{equation}
\label{int_{bf R} f(x) psi_{2^j}(x) dx = ....}
	\int_{\bf R} f(x) \, \psi_{2^j}(x) \, dx
   = \sum_{n=0}^\infty a_n \, 2^{-n \alpha} \, \widehat{\psi_{2^j}}(2^n)
	= a_j \, 2^{-j \alpha}.
\end{equation}
On the other hand,
\begin{equation}
	\int_{\bf R} \psi_{2^j}(x) \, dx = \widehat{\psi_{2^j}}(0)
		= \widehat{\psi}(0) = 0,
\end{equation}
so that
\begin{equation}
	\int_{\bf R} f(x) \, \psi_{2^j}(x) \, dx 
		= \int_{\bf R} (f(x) - f(0)) \, \psi_{2^j}(x) \, dx.
\end{equation}
Therefore
\begin{eqnarray}
	\Bigl|\int_{\bf R} f(x) \, \psi_{2^j}(x) \, dx \Bigr|
	& \le & \int_{\bf R} |f(x) - f(0)| \, |\psi_{2^j}(x)| \, dx	\\
  & \le & \|f\|_{\Lip \alpha} \int_{\bf R} |x|^\alpha \, |\psi_{2^j}(x)| \, dx
								\nonumber \\
  & = & \|f\|_{\Lip \alpha} \, 2^{-j \alpha} \, 
		\int_{\bf R} |x|^\alpha \, |\psi(x)| \, dx.	\nonumber
\end{eqnarray}
Combining this with (\ref{int_{bf R} f(x) psi_{2^j}(x) dx = ....}),
we obtain that
\begin{equation}
	|a_j| \le \|f\|_{\Lip \alpha} \, 
		\int_{\bf R} |x|^\alpha \, |\psi(x)| \, dx
\end{equation}
for all nonnegative integers $j$.  The integral on the right side
converges, because of the decay property of $\psi$.

	If $\alpha = 1$, then let us pass to the derivative and write
\begin{equation}
	f'(x) = \sum_{n=0}^\infty a_n \, i \, \exp{2^n \, i \, x}
\end{equation}
(where one should be careful about the meaning of $f'$ and of
this series).  This leads to
\begin{equation}
	\frac{1}{2\pi} \int_0^{2\pi} |f'(x)|^2 \, dx
		= \sum_{n=0}^\infty |a_n|^2.
\end{equation}
The main idea is that 
\begin{equation}
	\sum_{n=0}^\infty |a_n|^2 \le \|f\|_{\Lip 1}^2
\end{equation}
if $f$ is Lipschitz.  Conversely, if $\sum_{n=0}^\infty |a_n|^2 <
\infty$, then the derivative of $f$ exists in an $L^2$ sense, and
in fact one can show that $f'$ has ``vanishing mean oscillation''.

\section{Sums on general metric spaces}
\label{Sums on general metric spaces}
\setcounter{equation}{0}

	Let $(M, d(x,y))$ be a metric space.  For each integer $n$,
suppose that we have chosen a complex-valued Lipschitz function
$\beta_n(x)$ such that
\begin{equation}
	\sup_{x \in M} |\beta_n(x)| \le 1 \quad\hbox{and}\quad
		\|\beta\|_{\Lip } \le 2^n.
\end{equation}
Fix a real number $\alpha$, $0 < \alpha < 1$.

	Let $a_n$, $n \in {\bf Z}$ be a family (or doubly-infinite sequence)
of complex numbers which is bounded, and set 
\begin{equation}
	A = \sup_{n \in {\bf Z}} |a_n|.
\end{equation}
Consider 
\begin{equation}
	f(x) = \sum_{n \in {\bf Z}} a_n \, 2^{-n \alpha} \, \beta_n(x).
\end{equation}
The sum on the right side does not really converge in general, although
it would if we restricted ourselves to $n$ greater than any fixed number,
because of the bound on $\beta_n(x)$.  However, this sum does converge
``modulo constants'', in the sense that the sum in 
\begin{equation}
	f(x) - f(y) 
   = \sum_{n \in {\bf Z}} a_n \, 2^{-n \alpha} \, (\beta_n(x) - \beta_n(y)),
\end{equation}
converges absolutely for all $x$, $y$ in $M$.

	To see this, suppose that $k$ is any integer.  For $n \ge k$
we have that
\begin{equation}
	\sum_{n = k}^\infty |a_n| \, 2^{-n \alpha} \, |\beta_n(x)|
		\le A \, (1 - 2^{-\alpha})^{-1} \, 2^{-k \alpha},
\end{equation}
and similarly for $y$ instead of $x$.  For $n \le k-1$ we have that
\begin{eqnarray}
	\quad
   \sum_{n=-\infty}^{k-1} |a_n| \, 2^{-n \alpha} \, |\beta_n(x) - \beta_n(y)|
	& \le & A \sum_{n=-\infty}^{k-1} 2^{n (1-\alpha)} \, d(x,y)	\\
	& = & A \, 2^{(k-1)(1-\alpha)} \, (1 - 2^{-(1-\alpha)})^{-1} \, d(x,y).
								\nonumber
\end{eqnarray}
Thus
\begin{eqnarray}
\lefteqn{\sum_{n \in {\bf Z}} |a_n| \, 2^{-n \alpha} 
	   \, |\beta_n(x) - \beta_n(y)|}				\\
   & & \le A \, (1 - 2^{-\alpha})^{-1} \, 2^{-k \alpha}
     + A \, 2^{(k-1)(1-\alpha)} \, (1 - 2^{-(1-\alpha)})^{-1} \, d(x,y)
								\nonumber
\end{eqnarray}
for all $x, y \in M$ and $k \in {\bf Z}$.

\section{The Zygmund class}
\label{The Zygmund class}
\setcounter{equation}{0}

	Let $f(x)$ be a real or complex-valued function on the real
line.  We say that $f$ lies in the \emph{Zygmund class} $Z$ if $f$ is
continuous and there is a nonnegative real number $L$ such that
\begin{equation}
\label{def of Zygmund condition}
	|f(x+h) + f(x-h) - 2 \, f(x)| \le L \, |h|
\end{equation}
for all $x, y \in {\bf R}$.  In this case, the seminorm $\|f\|_Z$ is
defined to be the supremum of
\begin{equation}
	\frac{|f(x+h) + f(x-h) - 2 \, f(x)|}{|h|}
\end{equation}
over all $x, h \in {\bf R}$ with $h \ne 0$.  This is the same as the
smallest $L$ so that (\ref{def of Zygmund condition}) holds.  Clearly
$f$ is in the Zygmund class when $f$ is Lipschitz (of order $1$), with
$\|f\|_Z \le 2 \, \|f\|_{\Lip }$.

	Suppose that $\{a_n\}_{n=0}^\infty$ is a bounded sequence of
complex numbers, and consider the function $f(x)$ on ${\bf R}$
defined by
\begin{equation}
	f(x) = \sum_{n=0}^\infty a_n \, 2^{-n} \, \exp(2^n \, i \, x).
\end{equation}
Let us check that $f$ lies in the Zygmund class, with $\|f\|_Z$ bounded
in terms of
\begin{equation}
	A = \sup_{n \ge 0} |a_n|.
\end{equation}
Note that $f$ is continuous.

	Observe that
\begin{eqnarray}
\lefteqn{|\exp(i (u+v)) + \exp(i (u-v)) - 2 \exp (i \, u)|}	\\
	& & = |\exp (i \, v) + \exp (-i \, v) - 2|		\nonumber
\end{eqnarray}
for all real numbers $u$, $v$, and that
\begin{equation}
	\exp (i \, v) + \exp (-i \, v) - 2 
		= \int_0^v i (\exp(i \, t) - \exp (-i \, t)) \, dt
\end{equation}
when $v \ge 0$.  Since $|\exp (i \, t) - \exp(-i \, t)| \le 2 \, t$ for
$t \ge 0$, we obtain that
\begin{equation}
	|\exp (i \, v) + \exp (-i \, v) - 2| \le \int_0^v 2 \, t \, dt = v^2.
\end{equation}
Hence
\begin{equation}
\label{|exp(i(u+v))+exp(i(u-v))-2 exp(i u)| le v^2}
	|\exp(i (u+v)) + \exp(i (u-v)) - 2 \exp (i \, u)|
		\le v^2,
\end{equation}
and this works for all real numbers $u$, $v$, since there is no real
difference between $v \ge 0$ and $v \le 0$.

	Let $x$ and $h$ be real numbers, and let $m$ be a nonnegative
integer.  From (\ref{|exp(i(u+v))+exp(i(u-v))-2 exp(i u)| le v^2})
we get that
\begin{eqnarray}
\lefteqn{\Bigl| \sum_{n=0}^m a_n \, 2^{-n} \, 
 (\exp(2^n \, i (x+h)) + \exp(2^n \, i (x-h)) - 2 \exp(2^n \, i \, x)) \Bigr|}
									\\
   & & \qquad\qquad\qquad\qquad
   \le A \sum_{n=0}^m 2^{-n} \, 2^{2n} \, |h|^2 \le A \, 2^{m+1} \, |h|^2.
								\nonumber
\end{eqnarray}
If $|h| \ge 1/2$, then 
\begin{eqnarray}
\lefteqn{|f(x+h) + f(x-h) - 2 \, f(x)|}			\\
  & &  \le |f(x+h)| + |f(x-h)| + 2 \, |f(x)| \le 4 \, A \le 8 \, A \, |h|.
							\nonumber
\end{eqnarray}
If $|h| \le 1/2$, then choose a positive integer $m$ such that
$2^{-m-1} \le |h| \le 2^{-m}$.  We can write $f(x+h) + f(x-h) - 2 \, f(x)$
as
\begin{eqnarray}
\lefteqn{\qquad \sum_{n=0}^m a_n \, 2^{-n} \, 
 (\exp(2^n \, i (x+h)) + \exp(2^n \, i (x-h)) - 2 \exp(2^n \, i \, x)) }
									\\
  & & + \sum_{n=m+1}^\infty a_n \, 2^{-n} \, 
 (\exp(2^n \, i (x+h)) + \exp(2^n \, i (x-h)) - 2 \exp(2^n \, i \, x)).
								\nonumber
\end{eqnarray}
This leads to
\begin{eqnarray}
\lefteqn{|f(x+h) + f(x-h) - 2 \, f(x)|}			\\
  & &  \le |f(x+h)| + |f(x-h)| + 2 \, |f(x)|		\nonumber \\
  & &  \le A \, 2^{m+1} \, |h|^2 + 4 \, A \, 2^{-m}	\nonumber \\
  & &  \le A \cdot 2 \cdot |h| + 4 \cdot A \cdot 2 \cdot |h| = 10 \, A \, |h|.
							\nonumber
\end{eqnarray}
This shows that $f$ lies in the Zygmund class, with constant less than
or equal to $10 \, A$.

\section{Approximation operators, 1}
\label{Approximation operators, 1}
\setcounter{equation}{0}

	Let $(M, d(x,y))$ be a metric space.  Fix a real number
$\alpha$, $0 < \alpha < 1$, and let $f$ be a real-valued function
on $M$ which is Lipschitz of order $\alpha$.  For each positive
real number $L$, define $A_L(f)$ by
\begin{equation}
\label{def of A_L(f)}
	A_L(f)(x) = \inf \{f(w) + L \, d(x,w) : w \in M \}
\end{equation}
for all $x$ in $M$.

	For arbitrary $x$, $w$ in $M$ we have that
\begin{equation}
\label{f(w) ge f(x) - ||f||_{Lip alpha} d(x,w)^alpha}
	f(w) \ge f(x) - \|f\|_{\Lip \alpha} \, d(x,w)^\alpha.
\end{equation}
As a result,
\begin{equation}
	f(w) + L \, d(x,w) \ge f(x) 
\end{equation}
when $L \, d(x,w)^{1-\alpha} \ge \|f\|_{\Lip \alpha}$.  Thus we
can rewrite (\ref{def of A_L(f)}) as
\begin{eqnarray}
\label{def of A_L(f), 2}
\lefteqn{A_L(f)(x) = }					\\
	& & \inf \{f(w) + L \, d(x,w) : 
		w \in M, \ L \, d(x,w)^{1-\alpha} \le \|f\|_{\Lip \alpha} \},
							\nonumber
\end{eqnarray}
i.e., one gets the same infimum over this smaller range of $w$'s.  In
particular, the set of numbers whose infimum is under consideration is
bounded from below, so that the infimum is finite.

	Because we can take $w = x$ in the infimum, we automatically
have that
\begin{equation}
\label{A_L(f)(x) le f(x)}
	A_L(f)(x) \le f(x)
\end{equation}
for all $x$ in $M$.  In the other direction, (\ref{f(w) ge f(x) -
||f||_{Lip alpha} d(x,w)^alpha}) and (\ref{def of A_L(f), 2})
lead to
\begin{eqnarray}
	A_L(f)(x) & \ge & f(x) - \|f\|_{\Lip \alpha} \, 
	   \biggl(\frac{\|f\|_{\Lip \alpha}}{L} \biggr)^{\alpha/(1-\alpha)} \\
    & = & f(x) - \|f\|_{\Lip \alpha}^{1/(1-\alpha)} \, L^{-\alpha/(1-\alpha)}.
								\nonumber
\end{eqnarray}
We also have that $A_L(f)$ is $L$-Lipschitz on $M$, as in Lemma
\ref{lemma about inf of a family of L-Lipschitz functions}.

	Suppose that $h(x)$ is a real-valued function on $M$ which
is $L$-Lipschitz and satisfies $h(x) \le f(x)$ for all $x$ in $M$.
Then
\begin{equation}
	h(x) \le h(w) + L \, d(x,w) \le f(w) + L \, d(x,w)
\end{equation}
for all $x$, $w$ in $M$.  Hence 
\begin{equation}
	h(x) \le A_L(f)(x)
\end{equation}
for all $x$ in $M$.

	Similarly, one can consider
\begin{equation}
\label{def of B_L(f)}
	B_L(f)(x) = \sup \{f(w) - L \, d(x,w) : w \in M \},
\end{equation}
and show that
\begin{eqnarray}
\label{def of B_L(f), 2}
\lefteqn{B_L(f)(x) = }					\\
	& & \sup \{f(w) - L \, d(x,w) : 
		w \in M, \ L \, d(x,w)^{1-\alpha} \le \|f\|_{\Lip \alpha} \}.
							\nonumber
\end{eqnarray}
This makes it clear that the supremum is finite.  As before,
\begin{equation}
	f(x) \le B_L(f)(x) 
   \le f(x) + \|f\|_{\Lip \alpha}^{1/(1-\alpha)} \, L^{-\alpha/(1-\alpha)},
\end{equation}
and $B_L(f)$ is $L$-Lipschitz.  If $h(x)$ is a real-valued function
on $M$ which is $L$-Lipschitz and satisfies $f(x) \le h(x)$ for all $x$
in $M$, then
\begin{equation}
	B_L(f)(x) \le h(x)
\end{equation}
for all $x$ in $M$.

\section{Approximation operators, 2}
\label{Approximation operators, 2}
\setcounter{equation}{0}

	Let $(M, d(x,y))$ be a metric space, and let $\mu$ be a
positive Borel measure on $M$.  We shall assume that $\mu$ is a 
\emph{doubling} measure, which means that there is a positive
real number $C$ such that
\begin{equation}
	\mu(B(x,2r)) \le C \, \mu(B(x,r))
\end{equation}
for all $x$ in $M$ and positive real numbers $r$, and that the
$\mu$-measure of any open ball is positive and finite.

	Let $t$ be a positive real number.  Define a function
$p_t(x,y)$ on $M \times M$ by
\begin{eqnarray}
	p_t(x,y) & = & 1 - t^{-1} d(x,y) \quad\hbox{when } d(x,y) \le t	\\
		 & = & 0 	  \qquad\qquad\qquad\enspace
				       \hbox{when } d(x,y) > t,		
								\nonumber
\end{eqnarray}
and put
\begin{equation}
	\rho_t(x) = \int_M p_t(x,y) \, d\mu(y).
\end{equation}
This is positive for every $x$ in $M$, because of the properties of $\mu$.
Also put
\begin{equation}
	\phi_t(x,y) = \rho_t(x)^{-1} \, p_t(x,y),
\end{equation}
so that 
\begin{equation}
\label{int_M phi_t(x,y) dmu(y) = 1}
	\int_M \phi_t(x,y) \, d\mu(y) = 1
\end{equation}
for all $x$ in $M$ by construction.

	Fix a real number $\alpha$, $0 < \alpha \le 1$, and let $f$
be a complex-valued function on $M$ which is Lipschitz of order $\alpha$.
Define $P_t(f)$ on $M$ by
\begin{equation}
	P_t(f)(x) = \int_M \phi_t(x,y) \, f(y) \, d\mu(y).
\end{equation}
Because of (\ref{int_M phi_t(x,y) dmu(y) = 1}), 
\begin{equation}
	P_t(f)(x) - f(x) = \int_M \phi_t(x,y) \, (f(y) - f(x)) \, d\mu(y),
\end{equation}
and hence
\begin{eqnarray}
  |P_t(f)(x) - f(x)| & \le & \int_M \phi_t(x,y) \, |f(y) - f(x)| \, d\mu(y)
									\\
      & \le & \int_M \phi_t(x,y) \, \|f\|_{\Lip \alpha} \, t^\alpha \, d\mu(y)
		= \|f\|_{\Lip \alpha} \, t^\alpha.		\nonumber
\end{eqnarray}
In the second step we employ the fact that $\phi_t(x,y) = 0$ when
$d(x,y) \ge t$.

	Suppose that $x$ and $z$ are elements of $M$, and consider
\begin{equation}
	|P_t(f)(x) - P_t(f)(z)|.
\end{equation}
If $d(x,z) \ge t$, then
\begin{eqnarray}
\lefteqn{|P_t(f)(x) - P_t(f)(z)|}				\\
   & & \le |P_t(f)(x) - f(x)| + |f(x) - f(z)| + |P_t(f)(z) - f(z)|	
								\nonumber \\
   & & \le  \|f\|_{\Lip \alpha} (2 \, t^\alpha + d(x,z)^\alpha)
		\le 3 \, t^{\alpha - 1} \, \|f\|_{\Lip \alpha} \, d(x,z).
								\nonumber
\end{eqnarray}
Assume instead that $d(x,z) \le t$.  In this case we write $P_t(f)(x)
- P_t(f)(z)$ as
\begin{eqnarray}
\lefteqn{\int_M (\phi_t(x,y) - \phi_t(z,y)) \, f(y) \, d\mu(y)}		\\
    & & = \int_M (\phi_t(x,y) - \phi_t(z,y)) \, (f(y) - f(x)) \, d\mu(y),
								\nonumber
\end{eqnarray}
using (\ref{int_M phi_t(x,y) dmu(y) = 1}).  This yields
\begin{eqnarray}
\lefteqn{|P_t(f)(x) - P_t(f)(z)|}					\\
    & & \le \int_M |\phi_t(x,y) - \phi_t(z,y)| \, |f(y) - f(x)| \, d\mu(y)
								\nonumber \\
    & & \le (2 t)^\alpha \, \|f\|_{\Lip \alpha} 
	\int_{\overline{B}(x, 2t)} |\phi_t(x,y) - \phi_t(z,y)| \, d\mu(y),
								\nonumber
\end{eqnarray}
where the second step relies on the observation that $\phi_t(x,y) -
\phi_t(z,y)$ is supported, as a function of $y$, in the set
\begin{equation}
	\overline{B}(x,t) \cup \overline{B}(z,t) \subseteq \overline{B}(x, 2t).
\end{equation}

	Of course
\begin{eqnarray}
\lefteqn{\phi_t(x,y) - \phi_t(z,y)}					\\
    & & = (\rho_t(x)^{-1} - \rho_t(z)^{-1}) \, p_t(x,y)
		+ \rho_t(z)^{-1} \, (p_t(x,y) - p_t(z,y)).	\nonumber
\end{eqnarray}
Notice that
\begin{equation}
	|p_t(x,y) - p_t(z,y)| \le t^{-1} \, d(x,z)
\end{equation}
for all $y$ in $M$.  To see this, it is convenient to write $p_t(u,v)$
as $\lambda_t(d(u,v))$, where $\lambda_t(r)$ is defined for $r \ge 0$
by $\lambda_t(r) = 1 - t^{-1} \, r$ when $0 \le r \le t$, and
$\lambda_t(r) = 0$ when $r \ge t$.  It is easy to check that
$\lambda_t$ is $t^{-1}$-Lipschitz, and hence $\lambda_t(d(u,v))$ is
$t^{-1}$-Lipschitz on $M$ as a function of $u$ for each fixed $v$,
since $d(u,v)$ is $1$-Lipschitz as a function of $u$ for each fixed
$v$.  These computations and the doubling condition for $\mu$ permit
one to show that
\begin{equation}
	\int_{\overline{B}(x, 2t)} |\phi_t(x,y) - \phi_t(z,y)| \, d\mu(y)
		\le C_1 \, t^{-1} \, d(x,z)
\end{equation}
for some positive real number $C_1$ which does not depend on $x$, $z$,
or $t$.  (Exercise.)  Altogether, we obtain that
\begin{equation}
	\|P_t(f)\|_{\Lip 1} 
    \le \max(3, 2^\alpha \, C_1) \, t^{\alpha - 1} \, \|f\|_{\Lip \alpha}.
\end{equation}

\section{A kind of Calder\'on--Zygmund decomposition related to 
Lipschitz functions}
\label{A CZ decomposition related to Lipschitz functions}
\setcounter{equation}{0}

	Let $(M, d(x,y))$ be a metric space, and let $f$ be a
real-valued function on $M$.  Consider the associated maximal function
\begin{equation}
	N(f)(x) = \sup_{y \in M \atop y \ne x} \frac{|f(y) - f(x)|}{d(y,x)},
\end{equation}
where this supremum may be $+\infty$.

	Let $L$ be a positive real number, and put
\begin{equation}
	F_L = \{x \in M : N(f)(x) \le L \}.
\end{equation}
We shall assume for the rest of this section that
\begin{equation}
	F_L \ne \emptyset.
\end{equation}
As in Section \ref{Approximation operators, 1}, define $A_L(f)$
by
\begin{equation}
\label{def of A_L(f) repeated}
	A_L(f)(x) = \inf \{f(w) + L \, d(x,w) : w \in M \}.
\end{equation}
We shall address the finiteness of this infimum in a moment.
As before,
\begin{equation}
	A_L(f)(x) \le f(x)
\end{equation}
for all $x$ in $M$.

	If $u$ is any element of $F_L$, then
\begin{equation}
\label{|f(y) - f(u)| le L d(y,u), u in F_L}
	|f(y) - f(u)| \le L \, d(y,u)
\end{equation}
for all $y$ in $M$.  Let $x$ and $w$ be arbitrary points in $M$.
The preceding inequality implies that
\begin{equation}
	f(u) \le f(w) + L \, d(u,w),
\end{equation}
and hence
\begin{eqnarray}
	f(u) - L \, d(x,u) & \le & f(w) + L \, (d(u,w) - d(x,u))	\\
		& \le & f(w) + L \, d(x,w),			\nonumber
\end{eqnarray}
by the triangle inequality.  This yields
\begin{equation}
\label{f(u) - L d(x,u) le A_L(f)(x)}
	f(u) - L \, d(x,u) \le A_L(f)(x),
\end{equation}
which includes the finiteness of $A_L(f)(x)$.  If we take $x = u$,
then we get $f(u) \le A_L(f)(u)$, so that
\begin{equation}
	f(u) = A_L(f)(u) \quad\hbox{for all } u \in F_L.
\end{equation}
For $x \not\in F_L$, we obtain
\begin{equation}
	f(x) - 2 L \, d(x,u) \le A_L(f)(x)
\end{equation}
for all $u$ in $F_L$, by combining (\ref{f(u) - L d(x,u) le
A_L(f)(x)}) and (\ref{|f(y) - f(u)| le L d(y,u), u in F_L}) with $y =
x$.  In other words,
\begin{equation}
	f(x) - A_L(f)(x) \le 2 L \dist(x, F_L).
\end{equation}
Note that $A_L(f)$ is $L$-Lipschitz on $M$, by Lemma \ref{lemma about
inf of a family of L-Lipschitz functions}.

	In the same way, if
\begin{equation}
\label{def of B_L(f) repeated}
	B_L(f)(x) = \sup \{f(w) - L \, d(x,w) : w \in M \},
\end{equation}
then
\begin{equation}
	f(x) \le B_L(f)(x) \le f(x) + 2 L \dist(x,F_L)
\end{equation}
for all $x$ in $M$, and $B_L(f)$ is $L$-Lipschitz.


\end{document}